\documentclass[12pt,a4paper]{amsart}

\usepackage{amssymb}
\makeatletter
\renewcommand{\@begintheorem}[2]{
\rm \trivlist \item [\hskip \labelsep {\bf #2\ \ #1.}]
                                }
\makeatother

\makeatletter

\DeclareFontFamily{U}{cyr}{}
\DeclareFontShape{U}{cyr}{m}{n}{
  <5> wncyr5 <6> wncyr6 <7> wncyr7 <8> wncyr8 <9> wncyr9 <10->
wncyr10}{}
\DeclareMathAlphabet{\mathcyr}{U}{cyr}{m}{n}

\input cyracc.def

\newcommand{\ZZ}{{\bf Z}}
\newcommand{\QQ}{{\bf Q}}
\newcommand{\RR}{{\bf R}}
\newcommand{\CC}{{\bf C}}

\newcommand{\FF}{{\bf F}}
\newcommand{\HH}{{\bf H}}
\newcommand{\PP}{{\bf P}}
\newcommand{\AAA}{{\bf A}}

\newcommand{\cA}{{\mathcal A}}

\newcommand{\Km}{\mbox{Km}}

\newcommand{\rmd}{\mbox{d}}
\newcommand{\NS}{\mathop{\rm NS}}
\newcommand{\Sym}{\mathop{\rm Sym}}

\newcommand{\Jac}{\mathop{\rm Jac}}
\newcommand{\p}{\mathfrak p}

\newcommand{\bes}{\begin{equation}}
\newcommand{\ees}{\end{equation}}

\title{Two moduli spaces of abelian fourfolds with an automorphism of order five}
\author{Bert van Geemen}
\address{Dipartimento di Matematica, Universit\`a di Milano,
Via Saldini 50, I-20133 Milano, Italia}
\email{lambertus.vangeemen@unimi.it}

\author{Matthias Sch\"utt}
\address{Institut f\"ur Algebraische Geometrie, Leibniz Universit\"at
  Hannover, Welfengarten 1, 30167 Hannover, Germany}
\email{schuett@math.uni-hannover.de}

\begin{document}

\begin{abstract}
We find explicit projective models of a compact Shimura curve and of a 
(non-compact) surface which are the moduli spaces of principally polarised abelian fourfolds with an automorphism of order five. The surface has a $24$-nodal canonical model in $\PP^4$ which is the complete intersection of two $S_5$-invariant cubics. It is dominated by a Hilbert modular surface and we give a modular interpretation for this.  We also determine the L-series of these varieties as well as those of several modular covers of the Shimura curve.
\end{abstract}

\maketitle

Shimura varieties parametrise abelian varieties whose Mumford-Tate group is 
contained in a given algebraic group. Often this implies that the corresponding abelian varieties have a non-trivial endomorphism ring. In the simplest cases, Shimura varieties parametrise abelian varieties with a given automorphism group. In this paper we consider abelian fourfolds with an automorphism of order five.
The classical thetanulls map the moduli space of principally polarised abelian fourfolds with a level structure to projective space. 
As the automorphism of the abelian variety identifies various level structures
and the map is equivariant, the image of the Shimura variety must be in the fixed point locus in projective space of a specific projectivity of order five. 
In the cases we consider, this suffices to determine the Shimura varieties explicitly. 

We use the classical moduli space $A_{4,(2,4)}$ of principally polarised abelian fourfolds with a $(2,4)$-level structure. It has a nice equivariant map to a $\PP^{15}$ which is birational on its image.
It is also easy to find the action of projectivities corresponding to automorphisms on $\PP^{15}$ using the results from \cite{CG}. 
The equations for the image of the moduli space can be found from classical theta relations.
Using a computer, one then finds the equations for the intersection of the image of $A_{4,(2,4)}$ with certain eigenspaces of  specific projectivities of order five. These are then the equations for the images of the Shimura varieties, 
which in our cases turn out to be hypersurfaces in $\PP^2$ and $\PP^3$ respectively.
A similar method was already used in \cite{pmv} to find explicit models of Shimura
varieties parametrising abelian varieties with an automorphism of order three or four.
In those cases the projectivities were found by specialising to a product of elliptic curves, but this method does not work for automorphisms of order five.

We investigate the geometry of these varieties, in particular we determine the points which correspond to Jacobians of genus four curves, and we consider the number of vanishing theta\-nulls for the ppav's parametrised by these Shimura varieties.
We also determine the canonical model for the Shimura surface, it turns out to be the complete intersection of two $S_5$-invariant cubics in $\PP^4$.
The symmetric group $S_5$ has, up to tensoring by the sign representation, a unique irreducible 5-dimensional representation and this representation has a unique invariant cubic and a unique cubic which transforms with the sign representation. 
It is remarkable that the Shimura surface is the variety in $\PP^4$ defined by these two cubics, so it is `completely determined' by group theory.

The Shimura surface is dominated by a Hilbert modular surface. In fact, 
we show that the points on 
the Shimura surface correspond to abelian varieties which are isogenous 
to a self product of an abelian surface whose endomorphism algebra 
is an order in the field $\QQ(\sqrt{5})$.

We also consider the arithmetic of these Shimura varieties.
In particular, we investigate a tower of Shimura curves.
We show that the Jacobians of all curves in this tower
decompose completely into elliptic curves.
The L-series of these elliptic curves are explicitly related to classical and Hilbert modular forms.
A quotient of the Shimura surface is a K3 surface with a genus one fibration. 
We relate this fibration to the self product of an elliptic curve
from the tower of Shimura curves. 
This allows us to determine the zeta function of the Shimura surface.

\section{The moduli spaces $A_{g,(2,4)}$ and theta functions}

\subsection{General outline}\label{geno}
The Shimura varieties we consider are moduli spaces of 
principally polarised abelian varieties (ppav) with an automorphism, that is of triples
$(X,L,\phi)$ where $L$ is a line bundle defining a principal polarisation on the abelian variety $X$
and $\phi:X\rightarrow X$ is an automorphism such that $\phi^*L$ is algebraically equivalent to $L$.

Let $A_{g,n}$ be the moduli space of $g$-dimensional ppav's with  symplectic level 
$n$-structure 
$\alpha:X[n]\stackrel{\cong}{\rightarrow}(\ZZ/n\ZZ)^{2g}$.
The group $G=Sp(2g,\ZZ/n\ZZ)$ acts on $A_{g,n}$ by 
$h\cdot [(X,L,\alpha)]:= [(X,L,h\alpha)]$ and $A_{g,n}/G=A_g$.

An automorphism $\phi$ of $(X,L)$ acts on a level $n$-structure by $\alpha\mapsto \alpha\phi^{-1}$. Let $h\in Sp(2g,\ZZ/n\ZZ)$ be defined by
$h:=\alpha\phi^{-1}\alpha^{-1}$, then 
$h\cdot [(X,L,\alpha)]=[(X,L,\alpha\phi^{-1})]$.
On the other hand, $\phi$ being an automorphism,
the points $[(X,L,\alpha)]$ and $[(X,L,\alpha\phi^{-1})]$ in $A_{g,n}$ are the same. Thus $[(X,L,\alpha)]$ is a fixed point for the action of $h\in G$ on $A_{g,n}$. 

By considering deformations of $(X,L,\phi)$, one finds that all such deformations map to the fixed point set of $h$ in $A_{g,n}$. Conversely any fixed point of $h$
corresponds to an abelian variety with an automorphism.
If a fixed point lies in the same connected component of the fixed point set of $h$ in $A_{g,n}$ as $[(X,L,\phi)]$, then 
this abelian variety, with the automorphism corrresponding to $h$, is a deformation of $(X,L,\phi)$.

Let $\Theta_n:A_{g,n}\rightarrow\PP^N$ be
a $G$-equivariant embedding,
so for any $h\in G$ there is a linear map $M_h\in Aut(\PP^N)$ such that $\Theta_n\circ h=M_h\circ \Theta_n$. Then the image of the connected component
of the fixed point set of $h$ which contains $[(X,L,\alpha)]$ lies in the eigenspace
of $M_h$ which contains the point $\Theta([(X,L,\alpha)])$. 
In the rest of the paper, we will actually consider a slightly different level structure, the $(2,4)$-level structure, but we do follow this method.

\subsection{The $\Theta$-map}\label{Thetamap}
We consider the moduli space $A_{4,(2,4)}$ of ppav's of dimension four with a 
level two (symmetric) theta structure. 
Working for any $g$ and over $\CC$ from now on, this moduli space is a 
quotient of the Siegel upper half space $\HH_g$:
$$
A_{g,(2,4)}\,\cong\,\Gamma_g(2,4)\backslash \HH_g
$$
where $\Gamma_g(2)$ denotes the usual principal congruence subgroup of $\Gamma_g:=Sp(2g,\ZZ)$ of level $2$
\begin{eqnarray*}
\Gamma_g(2) & = & \ker(\Gamma_g\,\longrightarrow\,Sp(2g,\FF_2))\\
& = &
\left\{M=\begin{pmatrix}A&B\\C&D\end{pmatrix}\in Sp(2g,\ZZ):\;A\equiv D\equiv I,\;B\equiv C\equiv 0\;\mbox{mod}\;2\,\right\}
\end{eqnarray*}
and $\Gamma_g(2,4)$ is the (normal) subgroup of $\Gamma_g$
defined by
$$
\Gamma_g(2,4)\,:=\,
\left\{M=\begin{pmatrix}A&B\\C&D\end{pmatrix}\in \Gamma_g(2):\,
{\rm diag}A\,{}^t\!B \equiv {\rm diag }C\,{}^t\!D\equiv 0\,{\rm mod}\,4\,\right\}.
$$
We will write $A_\tau$ for ppav $\CC^g/\tau\ZZ^g+\ZZ^g$
with level structure defined by $\tau\in\HH_g$.

There are Galois covers:
$$
A_{g,4}\,\longrightarrow\,A_{g,(2,4)}\,\longrightarrow\,A_{g,2} \,\longrightarrow\,A_g,
$$
so a level two (symmetric) theta structure can be viewed as an equivalence class of level four structures. The cover $A_{g,(2,4)}\rightarrow A_{g,2}$ has group $\Gamma_g(2)/\Gamma_g(2,4)\cong \FF_2^{2g}$ and $A_{g,2} \rightarrow A_g$ has group $\Gamma_g/\Gamma_g(2)\cong Sp(2g,\FF_2)$, where $\FF_2$ is the field with two elements.

The advantage of this level structure is that, for any $g\geq 1$, there is a 
map
$$
\Theta:\,A_{g,(2,4)}\,\longrightarrow\,\PP^{2^g-1},\qquad 
[\tau]\,\longmapsto\,(\ldots:\Theta[\sigma](\tau):\ldots)
$$
where the thetanulls $\Theta[\sigma]$ for all $\sigma\in\{0,1\}^g$ are defined 
by the series 
$$
\Theta[\sigma](\tau)\,:=\,\theta_{\sigma/2\;0}(2\tau)\,=\,
\sum_{m\in\ZZ^g}\,e^{2\pi i {}^t(m+\sigma/2)\tau(m+\sigma/2)}\;.
$$
This map is a birational morphism onto its image.

In order to find the abelian varieties with automorphisms, 
we use
that the map $\Theta$ is equivariant for the action of 
$\Gamma_g$. 
So there is a  homomorphism
$$
\Gamma_g\,\longrightarrow\,PGL(2^g)\,:=\,GL(2^g)/\CC^\times,\qquad
M\,\longmapsto\, M_\Theta
$$
with kernel $\Gamma_g(2,4)$ such that
$$
\Theta(M\cdot\tau)\,=\,M_\Theta\Theta(\tau)\qquad(\forall \tau\in\HH_g).
$$
Thus the action of the $M_\Theta$ on $\Theta(A_{g,(2,4)})$  
corresponds to a change of the level structure.

\subsection{The action on projective space}\label{explicit}
The action of $\Gamma_g=Sp(2g,\ZZ)$ on $\PP^{2^g-1}$, which factors over $\Gamma_g/\Gamma_g(2,4)$,
can be described as follows.
The Galois group $\Gamma_g/\Gamma_g(2,4)$ of the cover 
$A_{g,(2,4)}\rightarrow A_g$ sits in an exact sequence (cf.\ \cite{DG}, 9.2):
{\renewcommand{\arraystretch}{1.7}
$$
\begin{array}{ccccccccc}
0&\longrightarrow&\Gamma_g(2)/\Gamma_g(2,4)&\longrightarrow&
\Gamma_g/\Gamma_g(2,4)&\longrightarrow&\Gamma_g/\Gamma_g(2)&
\longrightarrow&0\\
&&\cong\downarrow\phantom{\cong}
&&&&\cong\downarrow\phantom{\cong}&&\\
&&H_g/\mu_4 &&&&Sp(2g,\FF_2)&&
\end{array}
$$
}
The subgroup $\Gamma_g(2)/\Gamma_g(2,4)\cong\FF_2^{2g}$ acts on $\PP^{2^g-1}$ but this action does not lift to a linear action of this group on $\CC^{2^g}$.
Instead, a finite Heisenberg group $H_g$ acts on $\CC^{2^g}$ through its Schr\"odinger representation. The group $H_g$ is a central extension of 
$\FF_2^{2g}$ by the group of fourth roots of unity $\mu_4$. 
The group $\mu_4$ acts by scalar multiplication on $\CC^{2^g}$ hence $H_g/\mu_4=\FF_2^{2g}$ acts on $\PP^{2^g-1}$.

For $v\in \FF_2^{2g}$, $v\neq 0$, let $U_v\in GL(2^g,\CC)$ be the 
representation matrix  of an element in $H_g$ mapping to $v\in H_g/\mu_4$
such that $U_v^2=-I$.
We endow $\FF_2^{2g}$ with the standard symplectic form
$$
E\,:\FF_2^{2g}\times \FF_2^{2g}\longrightarrow \FF_2,
\qquad E((x,x'),(y,y'))\,:=\,{}^txy'+{}^tx'y
$$
for $x,x',y,y'\in\FF_2^g$. 
A point $v\in \FF_2^{2g}$  defines a transvection $t_v\in Sp(2g,\FF_2)$, which is an involution if $v\neq 0$,
by 
$$
t_v\,:\,V\,\longrightarrow\, V,\qquad t_v(w):=w+E(w,v)v.
$$
In \cite{CG} and \cite{DG} section 9.3 it is shown that, 
given such a transvection $t_v$, 
there is an $M\in \Gamma_g$ mapping to $t_v$ such that
$$
M_\Theta\,:\,\CC^{N+1}\,\longrightarrow\,\CC^{N+1},\qquad
M_\Theta\,:=\,\frac{1-i}{2}(U_v+I).
$$
If $M'\in \Gamma_g$ is another lift of $t_v$, then 
$M'_\Theta=U_wM_\Theta U_w^{-1}$
for some $w\in \FF_2^{2g}$. 
Thus it is easy to find explicitly the action of $\Gamma_g/\Gamma_g(2,4)$ on $\PP^{2^g-1}$.

\subsection{Products of ppav's}\label{products}
From the definition of the thetanulls $\Theta[\sigma]$ it is obvious that if 
the period matrix $\tau$ is in block form, with diagonal blocks $\tau_a,
\tau_b$ (with $\tau_c\in\HH_c$ and $a+b=g$) and zero entries elsewhere, then 
$$
\Theta[\sigma](\tau)=\Theta[\sigma_a](\tau_a) \Theta[\sigma_b](\tau_b)
$$
where if $\sigma=(s_1,\ldots,s_g)\in \{0,1\}^g$ then 
$\sigma_a=(s_1,\ldots,s_a)$, $\sigma_b=(s_{a+1},\ldots,s_g)$.

\subsection{Vanishing of thetanulls}
\label{vanth}
There are $2^{g-1}(2^g+1)$ even theta characteristics 
$[{}^\epsilon_{\epsilon'}]$ (that is with $\epsilon,{\epsilon'}\in \{0,1\}^g$
and $\epsilon\cdot{}^t\epsilon'=0$ mod $2$). For an even theta characteristic the corresponding thetanull $\theta[{}^\epsilon_{\epsilon'}](\tau)$, a holomorphic function on the Siegel space $\HH_g$,
does not vanish identically and it is a modular form of weight $1/2$ for a certain subgroup of $\Gamma_g$. 

The zero loci of the thetanulls have been studied extensively, especially in low genus. In particular, in points corresponding to products of ppav's (with product polarisation) and points in the boundary of the Satake compactification many thetanulls vanish (cf.\ \cite{pmv}, 3.6, 3.7). 
For example, let $A=A_\tau$ be a ppav which is the product of 
ppav's $B_1\times B_2$ of dimension $a$ and $b$, with period matrices $\tau_a,\tau_b$.
Then $\theta[{}^\epsilon_{\epsilon'}](\tau)
=\theta[{}^{\epsilon^{(a)}}_{{\epsilon'^{(a)}}}](\tau_a)
\theta[{}^{\epsilon^{(b)}}_{{\epsilon'^{(b)}}}](\tau_b)
$
and the $2^{a-1}(2^a-1)\cdot2^{b-1}(2^b-1)$ thetanulls 
with $\epsilon^{(a)}\cdot{}^t\epsilon'^{(a)}=1$ or $\epsilon^{(b)}\cdot{}^t\epsilon'^{(b)}=1$
vanish in $\tau\in \HH_g$. 

For example, $6\cdot 6=36$ even thetanulls vanish in a point
which is the product of two $2$-dimensional ppav's $B_1\times B_2$. 
In case $B_1$ is isomorphic to the product of two elliptic curves (with product polarisation), one thetanull of $B_1$ also vanishes, and this leads to the vanishing of $36+1\cdot 10=46$ thetanulls for the ppav $B_1\times B_2$, etc.

Any point in the boundary of (the Satake compactification of) $A_{g,(2,4)}$ is
the limit of a product of ppav's $B_1\times B_2$ where $B_2$ degenerates to a torus $(\CC^\times)^b$. In such a boundary point, besides the $2^{a-1}(2^a-1)\cdot2^{b-1}(2^b-1)$ thetanulls which vanish on the product, also 
the thetanulls $\theta[{}^{\epsilon^{(b)}}_{\epsilon'^{(b)}}]$
with at least one $\epsilon_i=1$ vanish. In particular, at least $100$ thetanulls vanish in a boundary point.

The following table lists the number of vanishing thetanulls for some cases of interest in $g=4$ (cf.\  \cite{pmv}, Lemma 3.8 for the cases $g=2,3$).

$$
\begin{array}{|@{\hspace{10pt}}l|@{\hspace{10pt}}c|}
\hline
\mbox{moduli point}& \mbox{\# vanishing thetanulls}
\\ \hline
B_1\times B_2, \quad B_i\; \mbox{abelian surfaces} & 36\\ \hline
(\CC^\times)^4 \quad \mbox{(boundary point)} & 120 \\ \hline
\end{array}
$$

\

In case the abelian fourfold is the Jacobian of a curve $C$
it is well-known that none of the thetanulls vanishes if $C$ is 
non-hyperelliptic and 
its image under the canonical map lies on a smooth quadric, 
exactly one vanishes if $C$ is non-hyperelliptic and its image 
under the canonical map lies on a rank three quadric, 
and exactly $10$ thetanulls vanish if $C$ is hyperelliptic.

\subsection{Thetanulls and quadrics}\label{thnq}
The classical formulas
$$
\theta[{}^\epsilon_{\epsilon'}](\tau)^2\,=\,
\sum_\sigma (-1)^{{}^t\sigma\epsilon'}\Theta[\sigma](\tau)\Theta[\sigma+\epsilon](\tau)
$$
show that 
there are quadrics $Q[{}^\epsilon_{\epsilon'}]$ in $\PP^{2^g-1}$ which pull back, along the map $\Theta$ to the zero locus of $\theta[{}^\epsilon_{\epsilon'}]^2$ in $\HH_g$, 
cf.\ \cite{pmv}, 3.5. 
The quadric $Q[{}^\epsilon_{\epsilon'}]$ is defined by the polynomial
$$
Q[{}^\epsilon_{\epsilon'}]\,:=\,Z(\,\sum_\sigma (-1)^{{}^t\sigma\epsilon'}
X_\sigma X_{\sigma+\epsilon}\,).
$$
Thus an abelian fourfold with a vanishing thetanull is mapped into one of these quadrics in $\PP^{15}$.

\subsection{Covers of $A_{g,(2,4)}$} \label{moco}

The relation between the $\theta[{}^\epsilon_{\epsilon'}]^2$ and the quadrics
$Q[{}^\epsilon_{\epsilon'}]$ recalled in section \ref{thnq} 
leads to a nice description of covers of the moduli space $A_{g,(2,4)}$. 
In fact, the map on $\HH_g$ given by both the $\Theta[\sigma]$ and the $\theta[{}^\epsilon_{\epsilon'}]$,
$$
\widetilde{\Theta}:\,\HH_g\,\longrightarrow\,\PP^M,
\qquad
\tau\,\longmapsto\,
(\ldots:\Theta[\sigma](\tau): \ldots :
\theta[{}^\epsilon_{\epsilon'}](\tau) :\ldots)
\qquad (M+1=2^g+2^{g-1}(2^g+1))
$$
factors over a normal subgroup $\Gamma_g(2,4,8)$ of $Sp(2g,\ZZ)$ 
(see \cite{vGN}). 
Thus the double cover of $\Theta(\HH_g)\subset\PP^{2^g-1}$ defined in $\PP^{2^g}\subset\PP^M$ by the equation
$T^2=Q[{}^\epsilon_{\epsilon'}]$ is the image (under a suitable projection) 
of $\widetilde{\Theta}(\HH_g)$. In particular, it is birational to a Siegel modular variety.
One can of course extend this result for any subset of the $Q[{}^\epsilon_{\epsilon'}]$.

\subsection{The equations for the moduli space}\label{eqmod}
The theta constants $\theta[{}^\epsilon_{\epsilon'}](\tau)$ in genus two satisfy 
certain quartic identities, for example the following two:
$$
\theta[{}^{00}_{00}]^4-\theta[{}^{00}_{10}]^4-
\theta[{}^{10}_{00}]^4-\theta[{}^{11}_{11}]^4\,=\,0,\qquad
\theta[{}^{01}_{00}]^4-\theta[{}^{10}_{10}]^4-
\theta[{}^{10}_{00}]^4+\theta[{}^{10}_{01}]^4\,=\,0.
$$
Replacing the fourth powers of the genus two theta constants with the following products of four genus four theta constants gives relation valid for all $\tau\in\HH_4$: 
$$
\theta[{}^{\epsilon}_{\epsilon'}]^4\;\longmapsto\;
\theta[{}^{00\epsilon}_{00\epsilon'}]\theta[{}^{00\epsilon}_{01\epsilon'}]
\theta[{}^{00\epsilon}_{10\epsilon'}]\theta[{}^{00\epsilon}_{11\epsilon'}]
$$
The two relations one obtains are of the form $r_1\pm r_2\pm r_3\pm r_4$, taking the product over the $2^3=8$ possible signs
$$
\prod_{\pm,\pm,\pm} (r_1\pm r_2\pm r_3\pm r_4)
$$
gives a relation which is a polynomial in $r_i^2$'s and thus is a polynomial in the squares of the theta constants in genus four. Now one uses the formula
from section \ref{thnq} to write each
$\theta[{}^\epsilon_{\epsilon'}]^2$ as a quadratic polynomial in the
$\Theta[\sigma]$
to obtain two equations, of degree 32, between the $ \Theta[\sigma]$'s.

\subsection{Schottky's modular form}\label{schottky}
The closure of the locus of the period matrices of Jacobians of curves in $\HH_4$
is the zero locus of
Schottky's modular form $J$. Explicitly, this modular form is defined by (cf.\  \cite{Ichr}):
$$
J\,:=\,2^4\sum\theta\left[{}^\epsilon_{\epsilon'}\right]^{16}\,-\,
\left(\sum\theta\left[{}^\epsilon_{\epsilon'}\right]^8\right)^2,
$$
where the sums are over the $136$ even theta characteristics.
Using the relation with the quadrics in $\PP^{15}$ for $g=4$ given in 
section \ref{thnq}, one finds a polynomial $F_J$, of degree $16$ in the $X_\sigma$'s, whose zero locus in $\PP^{15}$ 
intersected with the variety $\Theta(\HH_4)$ is the image of the zero locus of $J$, 
that is, it intersects the moduli space in the Jacobi locus.

\

\section{The four dimensional ppav's with order five automorphism}

\subsection{The ppav's $A_0$ and $B_0$}
Let $B_0:=J(C)$, be the Jacobian of the (unique) genus two curve $C$ with an automorphism of order five $\phi_C$:
$$
C:\quad y^2=x^5+1,\qquad 
\phi_{C}:\, (x,y)\,\longmapsto \,(\zeta x,y),
$$
where $\zeta$ is a primitive fifth root of unity.
Then $B_0$ is a principally polarized abelian surface. 
An p.p.\ abelian surface with an automorphism of order five must be simple
(if $E_1,E_2$ are elliptic curves, then $End(E_1\times E_2)\otimes\QQ$ cannot contain the field $\QQ(\zeta)$).
Thus it is the Jacobian of a curve with an automorphism of order five, 
so it is isomorphic to $B_0$.

Let $\phi_0\in End(B_0)$ be the automorphism induced by $\phi_C$, 
which thus preserves the principal polarization. 
The eigenvalues of $\phi_C$ on $H^{1,0}(B_0)\cong H^{1,0}(C)$ are 
$\zeta,\zeta^2$ since $\rmd x/y,x\rmd x/y$ are a basis of $H^{1,0}(C)$.

Let  $A_0:=B_0\times B_0$.
This is a p.p.\ abelian fourfold and it has 
automorphisms $\phi_k$:
$$
A_0:=B_0\times B_0,\qquad \phi_k:=\phi\times\phi^k\,\in\,End(A_0)
$$
which have eigenvalues $\zeta,\zeta^2,\zeta^k,\zeta^{2k}$
on $H^{1,0}(A_0)$. By deformation theory or \cite{BL}, section 9.6,
the deformation space of the pair $(A_0,\phi_k)$ has dimension 
$n_k=\dim (\Sym^2 H^{1,0}(A_0))^{\phi_k^*=1}$.
That is,
$n_k$ is the number of $1$'s in the set of $10$ products 
$\{\lambda_i\lambda_j\}_{1\leq i\leq j\leq k}$ 
where the $\lambda_i$ are the eigenvalues of $\phi_k^*$ on $H^{1,0}(A_0)$.
The deformation space is a Shimura variety in $A_4$, the moduli space of abelian fourfolds.

In case $k=1$, the pair $(A_0,\phi_1)$ is rigid. 
The cases $k=2,3$ are essentially the same, 
the deformations of $(A_0,\phi_2)$ are parametrised by a Shimura curve.
In case $k=4$, the eigenvalues are $\zeta,\zeta^2,\zeta^4=\overline{\zeta},
\zeta^8=\overline{\zeta}^2$, hence 
the deformations of $(A_0,\phi_4)$ are parametrised by a two dimensional 
Shimura variety, i.e.\ by a Shimura surface.

\subsection{The automorphisms $M_\Theta^{(k)}$} 
In case $g=2$, $A_{2,(2,4)}\cong\Theta(\HH_2)$ is the Zariski 
open subset of $\PP^3$ which is the complement of the $10$ quadrics $Q[{}^\epsilon_{\epsilon'}]$. 
The group $\Gamma_2/\Gamma_2(2,4)$ maps onto the group $\Gamma_2/\Gamma_2(2)\cong Sp(4,\FF_2)$ which is isomorphic to the symmetric group $S_6$, cf.\ \cite{CG}. Under this isomorphism, transpositions correspond to the transvections.
Thus it is easy to find an explicit element $M_\Theta\in PGL(4)$ 
of order five in the image of $\Gamma_2$. 
It has four one-dimensional eigenspaces in $\CC^4$, i.e.\ 
it has four fixed points in $\PP^3$ which do not lie in a quadric $Q[{}^\epsilon_{\epsilon'}]$.
Thus each of these points corresponds to the abelian surface $B_0$
with a choice of level structure.

We will also write $M_\Theta$ for a lift of $M_\Theta\in PGL(4)$ 
to $GL(4)$ which has  eigenvalues $\zeta,\ldots,\zeta^4$
and we denote by 
${\bf f}_1,\ldots,{\bf f}_4\in\CC^4$ the corresponding eigenvectors.

As the theta constants of the product $A_0=B_0\times B_0$ are products of the theta constants of $B_0$ (cf. section \ref{products}), we find that
any of the $4^2=16$ points ${\bf f}_i\otimes {\bf f}_j$ in $\PP(\CC^{16})=\PP(\CC^4\otimes\CC^4)$
corresponds to the fourfold $A_0$. 

As $\phi\in End(B_0)$ induces $M_\Theta$, the automorphism 
$\phi_k=\phi\times\phi^k$ induces 
$$
M^{(k)}_\Theta:=M_\Theta\otimes M^k_\Theta
\qquad (\in\,PGL(16)).
$$
Note that ${\bf f}_i\otimes{\bf f}_j$ is an eigenvector of 
$M^{(k)}_\Theta$ with eigenvalue $\zeta^{i+kj}$.

\subsection{The eigenspaces of $M^{(k)}_\Theta$}\label{eigm}
The point ${\bf f}_1\otimes{\bf f}_1$ in $\Theta(\HH_4)\subset \PP^{15}$ corresponds to $A_0$. 
The deformations of $(A_0,\phi_k)$ map to the 
eigenspace of $M^{(k)}_\Theta$ which contains this point,
which is the eigenspace with eigenvalue $\zeta^{1+k}$.

In case $k=2$, this eigenspace has basis 
$$
{\bf f_1}\otimes{\bf f_1},\quad {\bf f_2}\otimes{\bf f_3},\quad
{\bf f_4}\otimes{\bf f_2}.
$$ 
Thus the Shimura curve lies in the intersection 
of a $\PP^2$ with the variety $\Theta(\HH_4)$. 
Note that this curve will contain these three points, which all correspond to $A_0$ (with certain level structures).

In case $k=4$ this eigenspace  has basis 
$$
{\bf f_1}\otimes{\bf f_1},\quad
{\bf f_2}\otimes{\bf f_2},\quad
{\bf f_3}\otimes{\bf f_3},\quad
{\bf f_4}\otimes{\bf f_4}.
$$ 
Thus the Shimura surface lies in the intersection 
of a $\PP^3$ with the variety $\Theta(\HH_4)$. 

\subsection{The equations for the Shimura variety}\label{sheqs}
We are particularly lucky: 
in each case the Shimura variety is a hypersurface and is thus defined by just one equation. 
The equation was found by restricting two of the degree 32 polynomials on $\PP^{15}$ which vanish on $\Theta(\HH_4)$ as in section \ref{eqmod} to the eigenspace $\PP^n$, $n=2,3$. 
The greatest common divisor of these polynomials was irreducible in each case and is thus the defining equation of the Shimura variety. See sections \ref{vanshc} and \ref{eqshs}.

\

\section{Genus 4 curves and limit covers}
\label{s:limit}

\subsection{}
One important aspect of our investigations will be
the determination of the Jacobi locus inside our Shimura varieties.
Thus we are concerned with genus 4 curves
admitting an automorphism of order 5.
These come in one-dimensional families that we will describe below.
One of  these families is compact,
so if the curve degenerates, their Jacobians have limits which are still abelian fourfolds.
Here we sketch how to find these limit Jacobians through limit covers.

\subsection{Limit covers}

The general set-up of limit covers is as follows.
Consider an $n$-fold cover $C_\lambda$ of $\PP^1$ which depends on one parameter $\lambda$. 
Assume that $C_\lambda$ is locally around $(u,\lambda)=(0,0)$ given by
\[
C_\lambda\; : \; v^n \,=\, u^a\,(u-\lambda)^b
\]
for some $1\leq a,b< n$. 
We want to investigate the degeneration at $\lambda=0$. 
For simplicity, we restrict ourselves to the situation 
where $(a,n)=(b,n)=(a+b,n)=1$.
This implies that $n$ is odd;
upon degenerating, the genus of $C_\lambda$ drops by $(n-1)/2$.
Other cases can be treated along similar lines,
but the restriction is fairly natural, 
especially since we are aiming at the case $n=5$ anyway.

We will compensate the genus drop by introducing a curve of genus $(n-1)/2$ as a limit cover.
The construction goes as follows.
Consider the affine $(u,\lambda)$-plane $\AAA^2$. 
Then blow up an arbitrary point $P=(u_P,\lambda_P)$:
\[
 \tilde \AAA^2 \to \AAA^2.
\]
Denote the exceptional divisor by $E$. 
Then we consider the corresponding five-fold cover $\tilde X$ 
which is isomorphic to the total space $X$ of the family $C_\lambda$ 
outside the central fibre at $\lambda=\lambda_P$:
$$
\begin{array}{ccc}
 \tilde X & \dasharrow & X\\
\downarrow &&\downarrow\\
\tilde \AAA^2 & \longrightarrow & \AAA^2
\end{array}
$$
The limit cover involves the pre-image $\tilde E$ of the exceptional divisor $E$. 
The curve $\tilde E$ depends on whether $P$ lies on the lines $(u=0)$ or $(u=\lambda)$ in $\AAA^2$, 
since these lines constitute the local ramification locus.
We distinguish three cases:
\begin{enumerate}
 \item 
$P$ outside the lines $(u=0)$ and $(u=\lambda)$: 
Then $P$ has $n$ pre-images in $C_\lambda\times \AAA^1$, and so has $E$ in $\tilde X$, 
each isomorphic to $\PP^1$.
\item
$P\neq (0,0)$ on one of the lines $(u=0)$ or $(u=\lambda)$: 
Then $P$ is a ramification point. 
The pre-image of $E$ in $\tilde X$ is another $\PP^1$. 
The covering has degree $n$ and ramifies in the intersection points of $E$ with 
(the strict transforms of) the central fibre $(\lambda=\lambda_P)$ and the line.
\item
$P=(0,0)$: Again the pre-image of $E$ is a degree $n$ cover $\tilde E$ of $\PP^1$, 
but in this case there are three ramification points, 
namely the intersection points of $E$ with the central fibre 
(which in this case is the degenerate curve at $\lambda=0$) and the two lines.
\end{enumerate}

Our degenerate case concerns case (3).
In order to determine a local equation for $\tilde E$, 
we use the multiplicities of the branch points
which are given by the relevant exponents of the equation of $C$.
Presently a local equation for $\tilde E$ is given by
\[
 \tilde E:\;\; t^n \,=\, r^a\,(r-1)^b.
\]
By assumption, $\tilde E$ has genus $(n-1)/2$.
Thus we derived an interpretation of the degeneration of $\Jac(C_\lambda)$ at $\lambda=0$ as
$\Jac(C_0)\times \Jac(\tilde E)$.
The automorphism $v\mapsto \zeta_n v$ extends naturally to $\tilde E$ via $t\mapsto \zeta_n t$,
so we can even compute the induced action on the holomorphic 1-forms.

The above construction generalises to more complicated degenerations;
the details are left to the reader.
In the following, we concentrate on the case $n=5$.
We start with a family that will play an important role in \ref{ss:curve-Jac}.

\subsection{The de Jong--Noot family}
\label{ss:dJN}

In \cite{dJN}, the following family of genus $4$ curves was studied:
$$
C_\lambda:\;\; y^5=x(x-1)(x-\lambda).
$$
Let $\varrho$ denote a primitive sixth root of unity.
Then $C_\varrho$ is isomorphic to $y^5 = x^3-1$.
This implies that $\Jac(C_\varrho)$ has CM by $\QQ(\zeta_{15})$.
Starting from this Jacobian,
the authors show that the family has infinitely many CM points
(corresponding to simple abelian varieties).
This establishes a contradiction to  a conjecture of Coleman.

For completeness, we give a basis of $H^{1,0}$ 
consisting of eigenforms for the induced action of the automorphism 
$\phi: (x,y)\mapsto (x,\zeta y)$:
\[
\omega_1\,:=\,dx/y^4,\quad \eta_1\,:=\,x\omega_1,\quad
\omega_2\,:=\,y\omega_1,\quad \omega_3\,:=\,y^2\omega_1,
\]
where the index refers to the eigenvalues of $\phi^*$, so these are 
$\zeta, \zeta, \zeta^2, \zeta^3$. 
The canonical map $C\rightarrow \PP^3$ is thus given by 
$(x,y)\mapsto (X_0:\ldots:X_3)=(1:x:y:y^2)$
and the image of the curve lies in the rank three quadric 
$X_2^2=X_0X_3$. 
Thus the Jacobians of these curves have exactly one vanishing thetanull.

Now we study the degenerations at $\lambda=0, 1, \infty$ using limit covers.
At $\lambda=0$, the two 1-forms $dx/y^2, \;  xdx/y^4$ remain regular
with eigenvalues $\zeta^3, \zeta$.
Note that $C_0\cong C$.
The limit cover yields $\tilde E: t^5 = r(r-1)$ 
which is also isomorphic to $C$ with 
regular 1-forms $dt/r, tdt/r$ and eigenvalues $\zeta, \zeta^2$. 
Thus at $\lambda=0$, we find the degeneration $B_0\times B_0$ of $\Jac(C_\lambda)$
with automorphism $\phi_2$ (or $\phi_3$).
The argument at $\lambda=1$ is completely analogous and gives essentially the same result.

At $\lambda=\infty$, the degeneration of $C_\lambda$ is slightly different.
Namely $C_\infty\cong C$, 
but the regular 1-forms are $dx/y^3, dx/y^4$ with eigenvalues $\zeta^2, \zeta$.
On the other hand, the limit cover gives rise to
$\tilde E: t^5 = r^2(r-1)$
with regular 1-forms $dr/y^2, rdr/y^4$ and eigenvalues $\zeta^3, \zeta$.
Compared to the degenerations at $\lambda=0,1$,
degenerate curve and limit cover are interchanged,
but the product of the Jacobians still gives $B_0\times B_0$
with automorphism $\phi_2$ (or $\phi_3$).

\subsection{Cyclic five-fold covers of genus 4}

In genus two, there is a unique curve $C$ with an automorphism of order 5.
In this section we classify the cyclic five-fold covers of $\PP^1$ of genus $4$.
Essentially this amounts to a study of curves
\[
 v^5 = g(u)
\]
where $g(u)$ is a polynomial with four different root (including $\infty$).
Using symmetries, one can easily shrink the possibilities down to three families
up to isomorphisms:
the curves $C_\lambda$ from \ref{ss:dJN}
and the following two families:
\begin{eqnarray*}
 C'_\lambda: \;\; v^5 & = & u^4 (u^2+\lambda u+1),\\
C''_\lambda:\;\; v^5 & = & u^3 (u-1)^2 (u-\lambda).
\end{eqnarray*}
One easily computes the eigenvalues $(\zeta,\zeta^2,\zeta^3,\zeta^4)$
of the automorphism $\phi: v\mapsto \zeta v$ on the regular 1-forms.
This shows that two families $C'_\lambda, C''_\lambda$ are not isomorphic to $C_\lambda$.
In order to distinguish the two families,
we note that $C'_\lambda$ admits a hyperelliptic involution $(u,v)\mapsto (1/u, v/u^5)$.
In particular, the family $C'_\lambda$ can also be written canonically as 
\[
C'_\mu:\;\; y^2 = (x^5-1)(x^5-\mu).
\]
We point out that $C'_\lambda$ has a function of degree two, $f=u/v$, 
and the canonical map is $2:1$ onto the rational normal curve in $\PP^3$. 
On the other hand, it is easily verified using the regular 1-forms in \ref{ss:deg}
that the canonical map of $C''_\lambda$ is bijective on a Zariski open subset.
By the general theory, it follows that $C''_\lambda$ cannot be hyperelliptic.

It may be instructive to note that the curves $C_\lambda'$ and $C_\lambda''$ 
also appear prominently in the study of the Dwork pencil.
Indeed in \cite{candelas} it is observed that their $L$-functions are squares,
in agreement with our findings.

\subsection{Further degenerations}
\label{ss:deg}

It is an easy exercise to compute degenerations of the latest two families
of genus 4 curves using limit covers.
For $C'_\mu$, consider $\mu\to 0$.
From the 1-forms $x^idx/y \;(i=0,1,2,3)$, only $x^2dx/y, x^3dx/y$ remain regular on $C'_0$.
These have eigenvalues $\zeta^3, \zeta^4$.
As limit cover, we obtain
\[
 \tilde E: \;\; t^2 = r^5-1
\]
with induced eigenvalues $\zeta, \zeta^2$.
Again both genus two curves are isomorphic to $C$.
Thus $Jac(C'_\mu)$ degenerates as $B_0\times B_0$ at $\mu=0$ with automorphism $\phi_4$.

We turn to $C''_\lambda$ and consider again the degeneration $\lambda\to 0$.
The regular 1-forms of $C''_\lambda$ are spanned by the eigenforms
\[
 du/v,\; udu/v^2, u(u-1)du/v^3, u^2(u-1)du/v^4.
\]
Only the first two remain regular on $C''_0$. The limit cover gives 
\[
 \tilde E:\;\; t^5 = r^3(r-1)
\]
with regular 1-forms $r dr/t^3,\; r^2dr/t^4$. 
As above, $Jac(C'_\lambda)$ thus degenerates as $B_0\times B_0$ with automorphism $\phi_4$.

\section{The Shimura curve and its covers}

\subsection{The Shimura curve is a conic}
The Shimura curve is the intersection of an eigenspace
$\PP^2\subset\PP^{15}$ of $M^{(k)}_\Theta$ 
with $\Theta(\HH_4)\subset\PP^{15}$, whose equation, of degree two, was found as described in section \ref{sheqs}.
Thus it is a (smooth) conic in $\PP^2$ and so it has genus zero. 

\subsection{Vanishing of thetanulls on the Shimura curve}
\label{vanshc}
By computing the restrictions of the
quadrics $Q[{}^\epsilon_{\epsilon'}]$ in $\PP^{15}$ to the eigenspace $\PP^2$ we find the points on the curve where thetanulls vanish.

There is a unique quadric which vanishes on the $\PP^2$, thus the ppav's parametrized by the Shimura curve all have a vanishing thetanull.

The other $135=5\cdot 27$ quadrics restrict to $27$ conics in $\PP^2$, each of these conics is the restriction of 5 quadrics.
Actually $12$ of the conics are reducible and are the union of two distinct lines, the other $15$ are irreducible.
Each of these conics intersects the Shimura curve in $2\cdot 2=4$ points, counted with multiplicity. The union of the $27=12+15$ conics intersects the Shimura curve in only $12$ distinct points. 

The $12$ reducible conics intersect the Shimura curve in two points, one with multiplicity 1 and the other with multiplicity 3, in such a way that each of these 
twelve points determines a unique (reducible!) conic which has a triple zero there (so one of the two lines of the conic is tangent in that point) and there is a unique other point where this conic meets the curve. This sets up a pairing between the 12 points. 

Each of the $15$ irreducible conics meets the Shimura curve in $4$ distinct points, these four points form 2 pairs. Note that there are $6$ pairs of points and $\left({}^6_2\right)=15$ pairs of pairs of points, which is exactly the number of irreducible conics. Each pair of points lies in fact on $5$ irreducible conics.

In each of the 12 points exactly $36=1+5+5+5\cdot 5$ thetanulls vanish: one thetanull is zero on the Shimura curve, there are $5$ thetanulls which restrict to a reducible conic tangent to the point, there are $5$ thetanulls which restrict to a reducible conic not tangent to the point, and then there are 5 other conics which pass through the point. 

This implies that each of the $12$ points corresponds to a product of two abelian surfaces. As there are no points in the boundary of $A_4(2,4)$
where either exactly one or $36$ thetanulls vanish, the Shimura curve does not intersect the boundary, so it is a compact curve in $A_4(2,4)$.

\subsection{The Shimura curve lies in the Jacobi locus}
\label{ss:curve-Jac}

Using the polynomial $F_J$ from section \ref{schottky}, one finds
that the Shimura curve lies inside the closure of the moduli space of curves. 
As exactly one thetanull vanishes on the Shimura curve, the points of the
curve which are not in the boundary correspond to Jacobians of curves 
which are not hyperelliptic and whose canonical model lies on a cone (a singular quadric in $\PP^3$),
see \ref{vanth}.

In fact, we know from \ref{ss:dJN} 
 that the one parameter family of genus 4 curves $C_\lambda$ 
have a limit at $\lambda=0$
where its Jacobian specializes to $A_0=B_0^2$ 
and it has an automorphism of order five which specializes to $\phi_2$. 
Thus the Shimura curve is the moduli space of such curves (with a certain level structure).
In section \ref{ss:dJN} we observed that the image of the canonical map indeed lies on a cone in $\PP^3$.
We would like to point out 
that there are very few such (compact) Shimura curves of genus zero;
see \cite{J}.

\subsection{Parametrising the Shimura curve}
In our construction, the Shimura curve lies in a $\PP^2\subset\PP^{15}$
which is defined over the field of $20$-th roots of unity  $\QQ(i,\zeta)$. 
It is cut out by two polynomials of degree 16 with rational coefficients, 
so the curve, a conic, is defined over the field of $20$-th roots of unity. 
It has a rational point over that field. 
A convenient parametrisation puts the 6 pairs of points in the following position:
$$
\{0,\infty\},\quad\mbox{and}\quad  \{\zeta^k,\alpha\zeta^k\}_{k=0,\ldots,4}
\quad \mbox{with}\quad \zeta^5=1,\;\zeta\neq 1,\;\alpha:=\zeta^3+\zeta^2-1.
$$

\subsection{12 points on $\PP^1$}

One can ask how this 12-tuple of points sits in the moduli space
of 12 points in $\PP^1$.
We emphasize two aspects.
First 
one can appeal to the interpretation of the moduli space
in terms of isotrivial jacobian elliptic K3 surfaces
with j-invariant $0$.
That is, generically there are 12 singular fibers of type $II$ in Kodaira's notation
which determine the elliptic K3 surface over $\CC$.
In this setting, the 12-tuple is indeed very special
as it corresponds to the K3 surface of maximal Picard number $\rho=20$ 
and discriminant $-300$
studied in \cite{shioda}:
\begin{eqnarray}
\label{eq:300}
y^2 = x^3 + t^{11} - 11t^6-t.
\end{eqnarray}
We point out that the specified elliptic fibration has Mordell-Weil rank $18$;
 the given discriminant seems to be the smallest known in terms of absolute value
to allow for the maximal Mordell-Weil rank over $\CC$.

Surprisingly the Mordell-Weil group plays a crucial role for the second aspect.
Namely we were asked by D.~Testa
whether the 12 points
appear as discriminant of a rational elliptic surface,
that is, whether  the 12-tuple lies on the corresponding discriminant divisor in the moduli space.
In order to give an affirmative answer we note that
these rational elliptic surfaces correspond to integral sections of the K3 surface in \eqref{eq:300}.
Presently there are numerous such sections; they can be derived
 from the roots of the Mordell-Weil lattice $E_8$ of  the rational elliptic surface
 arising as quotient of \eqref{eq:300} induced by the involution $t\mapsto -1/t$ on the base curve $\PP^1$.
 Detailed computations of integral sections over $\QQ(\zeta_{15}, \sqrt[3]{10})$ are carried out in \cite[\S 5]{shioda}.

\subsection{Modular covers of the Shimura curve}
\label{ss:f}

As explained in section \ref{moco}, 
the covers of the Shimura curve given by equations of the type $t_i^2=Q_i$ 
where the $Q_i$ are the polynomials 
which vanish where the $Q[{}^\epsilon_{\epsilon'}]$ intersect the conic, 
are also modular curves.

In particular, the double cover of $\PP^1$ branched over the four points
$0,1,\infty,\alpha$ is also a Shimura curve.
The $j$-invariant of this elliptic curve is $16384/5$,
and it has a model over $\QQ$ in Weierstrass form with a rational point $(-1,1)$ of order 6:
\begin{eqnarray}
\label{eq:E}
E: \;\;\; y^2 = x(x^2+x-1).
\end{eqnarray}
The elliptic curve $E$ has conductor $20$.
Associated we find a cusp form of weight $2$ and level $20$:
\[
f = q - 2 q^3 - q^5 + 2 q^7 + q^9 + 2 q^{13} + 2 q^{15} - 6 q^{17} - 4 q^{19} - 4 q^{21} + 6 q^{23} + O(q^{25}).
\]

\subsection{Arithmetic of the covers of the Shimura curve}

The curves in the tower $t_i^2=Q_i \; (i=1,\hdots,r)$ admit a big automorphism group.
Immediately visible are $r$ copies of $\ZZ/2\ZZ$ which allow us to break 
up the Jacobian substantially.
Additional symmetries induce further involutions
which allow a further decomposition of the Jacobians.
In fact, we verified 
that the Jacobians always decompose into products of elliptic curves 
over the fixed extension $\QQ(i,\zeta)$.
The following table lists these elliptic curves
with their j-invariants,
the conductor of some preferred model over $\QQ$ or $\QQ(\sqrt 5)$ (such as \eqref{eq:E})
and the corresponding automorphic form.

\begin{table}[ht!]
$$
\begin{array}{ccc}
 \hline
\text{j-invariant} &
\text{conductor} &
\text{automorphic form}\\
\hline
2^{11}3^3/5 & 40 & g\\
2^{14}31^3/5^3 & 80 & f\otimes \chi_4\\
2^4 3^3 7^3/5^2 & 40 & g\\
\hline
2^7(25-11\sqrt 5)/5 & 40 & \mathfrak h\\
8(8903+3333\sqrt{5})/5 & 80 & \mathfrak h\otimes \chi_{4}\\
\hline
\end{array}
$$
\end{table}

Here $g=q + q^5 - 4q^7 - 3q^9 + 4q^{11} - 2q^{13} + 2q^{17} + 4q^{19} + 4q^{23} + O(q^{25})$
is the unique normalised newform of weight 2 and level $40$,
and $\chi_4$ denotes the quadratic Dirichlet character of conductor $4$.

The Hilbert modular form $\mathfrak h$ of parallel weight $2$ for $\QQ(\sqrt 5)$ is described as follows.
Consider the newspace inside $S_2(\Gamma_0(8\sqrt 5))$. 
A simple Magma computation reveals
that it has dimension 3
and is generated by the Hilbert modular form induced by $g$ and a pair of Galois-conjugate Hilbert modular forms
with Hecke eigenvalues in $\ZZ$.
Let $\psi$ denote the quadratic Hecke character of conductor $\sqrt 5$ and $\infty$-type $0$.
Twisting the two Hilbert modular forms by $\psi$ results in $\mathfrak h$ and its Galois-conjugate.
We list the first few Hecke eigenvalues $a_\p$ of $\mathfrak h$
according to the norm $\mathrm{N}(\p)$ of the prime ideal $\p$.

\begin{table}[ht!]
 $$
\begin{array}{cccccccccccc}
\hline
\mathrm{N}(\p) &&
 4 & 5 & 9 & 11 & 19 & 29 & 31 & 41 & 49 & 59\\
\hline
a_\p &&
0 & 0 & -2 & 4,-4 & 4,4 & -6, 10 & 8,0 & 10, -6 & 6 & -4, 12\\
\hline
\end{array}
$$
\end{table}

\subsection{}
\label{ss:a}

To prove the Hilbert modularity of the two elliptic curves over $\QQ(\sqrt 5)$,
one can proceed as follows.

The elliptic curves do not have 3-torsion over $\QQ(\sqrt 5)$, but
the reduction mod $3$ have by inspection of the trace $a_3=-2$.
Consider a quadratic twist, for instance over $\QQ(\sqrt\varepsilon)$ for $\varepsilon=(1+\sqrt 5)/2$,
such that the twisted curve $E'$ has trace $a_3=2$.
Then the argumentation from \cite{Dem} carries over to show that
\begin{itemize}
\item
the mod $3$ Galois representation $\bar\rho_{E',3}$ is irreducible and modular and
\item
the $3$-adic Galois representation is ordinary and absolutely irreducible so that 
\item
by \cite[Thm 5.1]{SW} the Hilbert modularity of $E'$ follows.
\end{itemize}
The level of the associated Hilbert modular form equals the conductor of $E'$ (which is $80$ in the above setting),
so we can find it by going through the full newspace of this level and comparing sufficiently many traces.
The Hilbert modular form associated to the original elliptic curve is obtained 
by twisting back
by the character corresponding to the above extension.

\subsection{}
Alternatively, one can
compute enough traces to apply the Faltings--Serre--Livn\'e method, see \cite{Livne}.
One subtlety here is the requirement that the Hilbert modular form $\mathfrak h$ 
(or its twist by the quadratic Hecke character of conductor $\sqrt 5$)
has eigenvalues in $\ZZ$ (which comes for free in \ref{ss:a}),
but in its current version Magma does not allow one to read off this property directly.
Once we know that the eigenvalues are in $\ZZ$,
it is easy to show that the traces are even and compute a test set of primes
so that equality of traces of Frobenius implies equality of $L$-functions.
For instance, the primes from \cite{DPS} will be more than sufficient for both aims
since they cover the case of good reduction outside ${2,3,\sqrt 5}$.

\section{The Shimura surface}

\subsection{The defining equations}\label{eqshs}
The Shimura  surface is the intersection 
of an eigenspace $\PP^3\subset\PP^{15}$ of $M^{(4)}_\Theta$
with $\Theta(\HH_4)\subset\PP^{15}$, whose equation, of degree six, was found as described in section \ref{sheqs}.
A study of the singular locus of this surface revealed an $S_5$-symmetry,
in fact, there are five `very singular' points, called the cusps, which are permuted by $S_5$ acting linearly on $\PP^3$.

In view of the $S_5$ symmetry, it was convenient
to introduce a $\PP^4$ which contains the $\PP^3$ as a hyperplane. The
equations of $S$ in this $\PP^4$ are: 
{\renewcommand{\arraystretch}{1.3}
$$
S:\qquad s_1:=x_1+\ldots+x_5\,=\,0,\qquad
s_2^3 + 10s_3^2 - 20s_2s_4\,=\,0,\qquad(\subset\PP^4)
$$
}
where $s_i$ is the $i$-th elementary symmetric function in the variables $x_1,\ldots,x_5$.

The singular points of $S$ are the 5 cusps in the orbit of $p_0$  and the 
$24$ points in the orbit of $q_0$ where
$$
p_0\,:=\,(-4:1:1:1:1),\qquad q_0\,:=\,(1:\zeta:\zeta^2:\zeta^3:\zeta^4),
$$
where $\zeta$ is a primitive 5-th root of unity (we used Magma to compute all singular points on $S$ over the algebraic closure of $\FF_{11}$). 
The tangent cone to $p_0$ is given by $xyz=0$ in $\PP^2$ whereas the tangent cone to $q_0$ is a smooth conic, so $q_0$ is node of $S$. As the singular locus has codimension 2, this also shows that $S$ is irreducible.

\subsection{The $22$ quadrics}

We consider the intersection of the 136 quadrics $Q[{}^\epsilon_{\epsilon'}]$ in $\PP^{15}$ 
with the eigenspace $\PP^3$
in which the Shimura surface $S$ lies.
Explicit computations show that these $136$ quadrics restrict to $22$ quadrics in $\PP^3$;
there are $6$ of these $22$ which are the restriction of $10$ quadrics in $\PP^{15}$
(we will call these type 10 quadrics), there are $15$ which are the restriction of $5$
(we call these type 5 quadrics) and there is one which is the restriction of a unique quadric (we call this quadric the invariant quadric as it is indeed $S_5$-invariant).
This accounts for the restriction of all the $136$ quadrics:
$$
136\,=\,6\cdot 10\,+\,15\cdot5\,+\,1.
$$
These 22 quadrics in $\PP^3$ are all defined over $\QQ$ in the coordinates on $\PP^4$ given above.

\subsection{Vanishing of thetanulls in the singular points}

The surface $S$ has one $S_5$-orbit of 5 singular points.
Each of these singular points
lies on $18$ of the $22$ quadrics, only the invariant quadric and three of the $15$  type $5$ quadrics do not pass through it. In particular, 
$120$ thetanulls vanish in such a point and thus it corresponds to a maximal degeneration of an abelian fourfold (denoted by $(\CC^\times)^4$ in the table in section \ref{vanth}).
We will call these 5 singular points of $S$ the cusps from now on.
The $6$ type 10 quadrics are the only ones from the $22$ to pass through all five cusps.

The other $24$ singular points, the $S_5$-orbit of $q_0$, 
are not defined over $\QQ$ but over $\QQ(\zeta)$. 
 The full orbit of 24 singular points
consists of $4\cdot 6$ Galois orbits.
As the $22$ quadrics are defined over $\QQ$ this implies that if 
one of them contains one of these singular points, 
then it contains the Galois orbit of that point
which consists of $4$ points. 

In $q_0$ exactly  $7$ of the $22$ quadrics vanish: 
one of type $10$, $5$ of type $5$ and the invariant quadric.
Thus $36$ thetanulls vanish in $q_0$ and 
therefore $q_0$ corresponds to a product of two abelian surfaces, 
like $A_0=B_0\times B_0$ (and we will see that $A_0$ produces these singular points). 
We will call the $24$ points in the orbit of $q_0$ the reducible points.
The vanishing of a quadric on a Galois orbit gives a bijection 
between the six type 10 quadrics and the six Galois orbits on the 
reducible points.

\subsection{Vanishing of thetanulls on the Shimura surface}
The type 10 quadrics are smooth quadrics in $\PP^3$, 
hence they are isomorphic to $\PP^1\times\PP^1$ and we can parametrise them.
In this way the intersection of a type 10 quadric with the Shimura surface,
of degree 6 in $\PP^3$, was determined explicitly. 
It is a bidegree $(6,6)$ curve on $\PP^1\times\PP^1$. 
Explicit computations show that this curve has two irreducible components over $\QQ(\sqrt 5)$, 
each of which has multiplicity 2, and the classes of the curves are $(1,2),(2,1)$. 
In particular these are rational curves (twisted cubics in $\PP^3$) 
which meet in $5=1\cdot1+2\cdot 2$ points, 
these are just the cusps of the Shimura surface.

Similarly, the type 5 quadrics are smooth. 
Their intersection with the Shimura surface is a $4$-nodal curve of bidegree $(3,3)$ with multiplicity two;
this curve is irreducibe and rational.

\subsection{Jacobians in the Shimura surface}

To find the points of $S$ which correspond to (limits of) Jacobians, 
we used the polynomial $F_J$ from section \ref{schottky}.

We considered first the intersection of the degree $16$ surface defined by $F_J=0$ 
with the type 10 quadrics.
It was easy to check (with MAGMA) that $F_J$ cuts a type 10 quadric
in three curves, two of these curves are the intersection of the quadric with $S$, and these have again multiplicity two, and there is a residual curve of degree 10 which does not lie in $S$. 
These two curves thus parametrize Jacobians with 10 vanishing theta nulls (the general point of such a  curve lies in only one of the 22 quadrics). Therefore these must be hyperelliptic Jacobians.
In fact, we already know from \ref{ss:deg}
that  the Jacobians of the hyperelliptic genus 4 curves $C'_\lambda$ 
specialise to $A_0=B_0^2$ as  
$\lambda\rightarrow 0$
 and that they
have an automorphism of the right type ($=(1,1,1,1)$) which specializes to $\phi_4$, so the two curves in the triple intersection of a type 10 quadric, $S$ and $F_J=0$ 
parametrise such hyperelliptic curves.

The intersection of $S$ with $F_J=0$ is a curve of degree $6\cdot 16=96$.
We have already seen that each of the six type 10 quadrics contains two curves of degree three of $S\cap (F_J=0)$, each with multiplicity two, thus there remains a curve of degree $96-6\cdot (3+3)\cdot 2=24$.
We verified, using MAGMA again, that this curve meets the type 5 quadrics only in a finite set of points. 
Thus the general point of this curve corresponds
to a non-hyperelliptic Jacobian, as it 
does not have a vanishing theta null. On the other hand, we know from \ref{ss:deg}
that the
Jacobians of the non-hyperelliptic curves $C''_\lambda$
specialize to $A_0=B_0^2$ and have an automorphism which specializes to $\phi_4$, so these are the Jacobians parametrised by the degree $24$ curve in the Shimura surface $S$.

\subsection{The canonical model of the Shimura surface}
The Shimura surface $S$ has degree $6$ in $\PP^3$, so the canonical system is cut out by quadrics. 
These quadrics must pass through the five cusps in order to define regular $2$-forms on the desingularisation of $S$.

A computation shows that the vector space $V$ of quadrics vanishing in the cusps has dimension $5$. 
It is spanned by the six type 10 quadrics, which are permuted under the action of $S_5$. The representation of $S_5$ on $V$ is irreducible 
and, up to tensoring with the alternating representation of $S_5$, it is the unique $5$-dimensional irreducible representation of $S_5$. 
The representation $V$ is self dual and thus the representation on the cubics on $V$ is the representation $S^3V$.
This representation has one invariant and one alternating cubic.
It turns out that the image of $S$ under the canonical map is the surface in $\PP V\cong \PP^4$ defined by these two cubics.  Explicitly, we identify $\PP V$ with a hyperplane in $\PP^5$. The action of $S_5$ is as a subgroup of $S_6$, which acts by permutation of the coordinates, but the embedding $S_5\hookrightarrow S_6$, corresponding to the action on the six type $10$ quadrics, is not the standard embedding, but is given by
$$
S_5\,\hookrightarrow \,S_6\qquad
(1\,2)\,\longmapsto\,(1\,4)(2\,3)(5\,6),\qquad
(5\,4\,3\,2\,1)\,\longmapsto\,(2\,6\,5\,4\,3).
$$
The equations defining the image of $S$ are then:
$$
z_1+z_2+z_3+z_4+z_5+z_6,\qquad
z_1^3+z_2^3+z_3^3+z_4^3+z_5^3+z_6^3,
$$
and the following (alternating for $S_5$) cubic:
\begin{eqnarray*}
z_1z_2z_3 - z_1z_2z_4 - z_1z_2z_5 + z_1z_2z_6 + z_1z_3z_4 - z_1z_3z_5 - z_1z_3z_6 + 
z_1z_4z_5 - z_1z_4z_6 + z_1z_5z_6 +
\\
- z_2z_3z_4 + z_2z_3z_5 
- z_2z_3z_6 + z_2z_4z_5 +
z_2z_4z_6 - z_2z_5z_6 - z_3z_4z_5 + z_3z_4z_6 + z_3z_5z_6 - z_4z_5z_6.
\end{eqnarray*}

Thus the image of the Shimura surface is the complete intersection of two cubics in a $\PP^4$. A computation shows that it has only 24 singular points, 
the five cusps are resolved into five cycles of three lines. 
In particular, it is a canonically embedded surface.
A computer computation suggests that the map from $S$ to its image in $\PP V$ is a birational morphism.

\section{The arithmetic of the Shimura surface}

\subsection{Picard number}

Our first aim is to compute the Picard number of the Shimura surface over $\CC$.
Throughout we work with the  minimal desingularization $\tilde S$ of the sextic model in $\PP^3$.
From the canonical model as a (3,3) surface in $\PP^4$,
one easily finds the following invariants:
\[
 q(\tilde S)=0,\;\; p_g(\tilde S)=5,\;\; K_{\tilde S}^2=9,\;\; e(\tilde S)=63,\;\; h^{1,1}(\tilde S)=51.
\]
By Lefschetz' theorem, one has $\rho(\tilde S)\leq h^{1,1}(\tilde S)=51$.

\subsection{Lemma}
{\sl The complex Shimura surface $\tilde S$ has Picard number $\rho(\tilde S)=46$.}
\label{ss:lem}

\medskip
 
The proof of the lemma proceeds in three steps \ref{ss:1st}--\ref{ss:3rd}.
First we exhibit a natural collection of 46 independent divisors on $\tilde S$,
so that $\rho(\tilde S)\geq 46$
Then we use the $S_5$ action to show that $\rho(\tilde S)>46$ would imply $\rho(\tilde S)=51$.
Moreover, we deduce that $\tilde S$ would be modular by a result of Livn\'e.
Then we establish a contradiction by counting points over small finite fields
and applying Lefschetz' fixed point formula.

\subsection{Independent divisors on $\tilde S$}
\label{ss:1st}

It suffices to consider a fairly simple collection of curves on $\tilde S$:
\begin{enumerate}
 \item 
the exceptional curves over the 24 nodes; these are $(-2)$ curves;
\item
the components of the exceptional divisors $E_i\; (i=1,\hdots,5)$ over the 5 cusps;
these are $(-3)$-curves forming a triangle;
\item
the two components over $\QQ(\sqrt 5)$ of the 6 type 10 quadrics;
these are $(-3)$ curves as well.
\end{enumerate}

The self-intersection numbers can be computed with the adjunction formula as follows.
The given model of $\tilde S$ as desingularisation of a singular sextic in $\PP^3$
has canonical divisor
\[
 K_{\tilde S} = 2H - (E_1+\hdots + E_5).
\]
For $K_{\tilde S}^2$ to be $9$, we require $E_i^2=-3$
which implies that each component is a $(-3)$ curve.
Meanwhile the two irreducible components $C_1, C_2$ of a type 10 quadric $Q$
meet exactly one component of each 
exceptional divisor over a cusp, each transversally in a single point.
Here we have $Q=2(C_1+C_2)$, so $C_1.H=Q.H/4=3$ by symmetry.
Thus $C_i.K_{\tilde S}=1$, and since $C_i$ is smooth and rational, $C_i^2=-3$ as claimed.

More precisely, each component of an exceptional divisor over a cusp
is met by two different pairs of conjugate twisted cubics,
the intersection being transversal in four distinct points.
In particular, this shows that all 24 twisted cubics 
become disjoint on the desingularisation.
This information allows us to write down the intersection matrix of the 51 given curves.
A machine calculation returns that the matrix has rank $46$,
so these curves generate a sublattice $N$ of $\NS(\tilde S)$ of rank 46.
We give two preferred basis of $N$:
\begin{enumerate}
 \item 
\label{eq:N}
Omit five twisted cubics that are pairwise non-conjugate.
The resulting Gram matrix has determinant $-2^{45}3$.
In particular, since the twisted cubics are interchanged by Gal$(\QQ(\sqrt 5)/\QQ)$,
this shows that as a Galois module 
\begin{eqnarray*}
 N=\QQ(-1)^{10}\oplus(\QQ(\sqrt 5))(-1)^6\oplus(\QQ(\zeta))(-1)^6.
\end{eqnarray*}
\item
Omit the exceptional curve over one node from each Galois orbit, 
but add the hyperplane section $H$.
The resulting Gram matrix has determinant $-2^{33}3$.
\end{enumerate}

\subsection{Transcendental lattice}

The transcendental  lattice $T(\tilde S)$ is defined as the orthogonal complement of $\NS(\tilde S)$ in $H^2(\tilde S,\ZZ)$
with respect to cup-product:
\[
 T(\tilde S) = \NS(\tilde S)^\bot \subset H^2(\tilde S,\ZZ).
\]
It comes equipped with a Hodge structure
and with a compatible system of Galois representations over $\QQ$.
Recall that the six type 10 quadrics (with one linear relation) 
represent the regular 2-forms on $\tilde S$.
Since $S_5$ acts as an irreducible representation on them,
this group action splits up $T(\tilde S)$ 
as a Hodge structure  or as a Galois representation over $\QQ$
\[
 T(\tilde S) = V^5.
\]
Here $V$ has Hodge-Tate weights $(2,0)$ and $(0,2)$ plus possibly $(1,1)$.
Since we already know that $\rho(\tilde S)\geq 46$,
we deduce that $V$ is either 2 or 3-dimensional and $\rho(\tilde S)=46$ or $51$.

\subsection{}
\label{ss:3rd}

In the present situation, 
the property  $\rho(\tilde S)=46$ could be deduced 
from the geometric argumentation in Section \ref{s:K3}.
For completeness, we give an ad-hoc proof that does not rely on any geometric insight.

Assume that $V$ has dimension $2$.
Then $V$ gives rise to a 2-dimensional motive over $\QQ$ 
endowed with a non-degenerate intersection form.
By a theorem of Livn\'e \cite{L} (or generally nowadays by Serre's conjecture),
$V$ would be modular.
The associated newform would have weight $3$ and rational Fourier coefficients $a_n$
such that
\[
 \mbox{trace Frob}_p^*|_V = a_p
\]
for almost all primes $p$.
By \cite[Lemma 1.8]{L}, this implies CM by an imaginary quadratic field $K$
whose class group is only 2-torsion.
The field $K$ can be recovered from the newform
by the roots $\alpha, \bar\alpha$ of the local Euler factors 
(which are the eigenvalues of Frobenius on cohomology).
Namely we always have $\alpha\in K$.

These ideas can be used to establish a contradiction to our assumption as follows:

By the Lefschetz fixed point formula, we have at a prime $p$ of good reduction
\begin{eqnarray}\label{eq:Lef}
 5\, \mbox{trace Frob}_p^*|_V = \# \tilde S(\FF_p) - 1 - p\, \mbox{trace Frob}_p^*|_{\NS(\tilde S)} - p^2.
\end{eqnarray}
Here $N$ gives $\NS(\tilde S)$ up to a 5-dimensional piece,
so we know $\mbox{trace Frob}_p^*|_{\NS(\tilde S)}$ up to a summand $h\in\{-5,\hdots,5\}$.
But then the divisibility condition on the LHS in \eqref{eq:Lef}
leaves only two or three possibilities for $h$.
In the present situation we find the following possibilities:
\begin{table}[ht!]
 $$
\begin{array}{cccc}
\hline
 p & h & \mbox{trace Frob}_p^*|_V & \alpha\\
\hline
3 & 5 & -2 & -1\pm 2\sqrt{-2}\\
& 0 & 1 & (1\pm\sqrt{-35})/2\\
& -5 & 4 & 2\pm\sqrt{-5}\\
\hline
7 & 5 & -10 & 5+2\sqrt{-6}\\
& 0 & -3 & (-3\pm\sqrt{-11\cdot 17})/2\\
& -5 & 4 & 2\pm 3\sqrt{-5}\\
\hline
13 & 5 & -22 & -11\pm 4\sqrt{-3}\\
& 0 & -9 & (-9\pm\sqrt{-5\cdot 7\cdot 17})/2\\
& -5 & 4 & 2\pm\sqrt{-3\cdot 5\cdot 11}\\
\hline
\end{array}
$$
\end{table}

Since there is no overlap between the possible fields $\QQ(\alpha)$
between $p=7$ and $p=13$,
we deduce that $V$ cannot be 2-dimensional.
Hence $\dim(V)=3$ and $\rho(\tilde S)=46$.
This completes the proof of Lemma \ref{ss:lem}.
\qed

\subsection{Remarks}

\begin{enumerate}
 \item[(a)]
If one prefers not to work with $T=V^5$,
one can also appeal to  quotient surfaces of $\tilde S$ such as the K3 surface $X$ in Section \ref{s:K3}.
$\tilde S$ admits also a free automorphism of order $3$. 
Both quotients have $p_g=1$, so exactly one copy of $V$ descends from $\tilde S$ to the quotients.
Thus one can argue directly with $V$ on the quotient to prove that it cannot have dimension 2.
\item[(b)]
Instead of appealing to modularity in \ref{ss:3rd},
one can also use the Tate conjecture for the reductions of $\tilde S$ (or $X$) mod $p$ 
and its refinement, the Artin-Tate conjecture
(without actually relying on its validity).
Essentially this amounts to checking the same information as above, see \cite{S-NS}
for general arguments towards singular K3 surfaces
or \cite{RS} for a similar example of a surface of general type which is not defined over $\QQ$.
\end{enumerate}

\subsection{Zeta function}

For $h=5$,
the Frobenius eigenvalues $\alpha$ in the table above equal 
the squares of the roots of the newform $f$ from \ref{ss:f}.
That is, if $f=\sum a_nq^n$, 
then the 3-dimensional Galois representation has eigenvalues $p, \alpha, \bar\alpha$ as above
such that
\[
 \mbox{trace Frob}_p^*|_V = p + \alpha + \bar\alpha = a_p^2 - p.
\]
Counting points over $\FF_p$ and $\FF_p^2$, one can verify this  for many further primes.
This suggests the following:

\subsection{Theorem}
\label{ss:thm}
{\sl
The following two L-functions coincide:
\[
 L(\mbox{Sym}^2 H^1(E),s) = L(V,s).
\]
}

We will give a geometric proof of the theorem in Section \ref{s:K3}
using the theory of K3 surfaces.
For the zeta function of $\tilde S$, we also have to take into account the algebraic cohomology pieces.
Since $\NS(\tilde S)\otimes\QQ = N\otimes\QQ$, these are dealt with by \ref{ss:1st} \eqref{eq:N}.
\subsection{Corollary}
{\sl 
The zeta function of $\tilde S$ is given as
\[
 \zeta(\tilde S,s) = \zeta(s) \zeta(s-1)^{10} \zeta_{\QQ(\sqrt 5)}(s-1)^6 
\zeta_{\QQ(\zeta)}(s-1)^6 L(\mbox{Sym}^2 H^1(E),s)^5 \zeta(s-2).
\]}

\section{Relation with K3 surfaces}
\label{s:K3}

This section is devoted to a proof of Theorem \ref{ss:thm}.
We will proceed as follows:
first we exhibit a quotient $X$ of $\tilde S$ with $p_g=1$.
Using elliptic fibrations, we  relate $X$ in several steps to a K3 surface of Kummer type,
namely for the elliptic curve $E$ from \ref{ss:f}
and its 6-isogenous partner $E'$.
Our approach relies on Shioda-Inose structures.
Each intermediate step will be exhibited over $\QQ$.

\subsection{K3 quotient}

On $\tilde S$, consider the automorphism $\imath$ of order 5 permuting coordinates cyclically.
Taking eigencoordinates and cubic invariants,
one finds the quotient $X=\tilde S/\imath$ as a sextic in $\PP^3$ given by the polynomial
\begin{eqnarray*}
4 z_1^3 z_2 z_4^2-3 z_1^3 z_4^3-z_1^2 z_2 z_3 z_4^2+4 z_1^2 z_3 z_4^3-2 z_1 z_2^3 z_3 z_4+4 z_1 z_2^2 z_3^3-
\\
z_1 z_2^2 z_3^2 z_4-2 z_1 z_2 z_3^3 z_4-2 z_1 z_2 z_3 z_4^3-3 z_2^3 z_3^3+4 z_2^3 z_3^2 z_4-2 z_1^3 z_2 z_3 z_4.
\end{eqnarray*}
The sextic $X\subset\PP^3$ admits an automorphism $\jmath$ of order 4, permuting coordinates through $(1243)$.
$X$ is singular along the 4 lines which form the orbit of $(z_3=z_4=0)$ under $\jmath$.
Moreover there is a triple point at $[1,1,1,1]$.
Any regular 2-form has to vanish along these lines and in the triple point,
so the unique regular 2-form is given by the quadric $z_1z_4-z_2z_3$ up to scalars.
With a good portion of work,
one can show that the minimal desingularization of $X$ can be blown down to a K3 surface.
Here we will not pursue this argument,
but rather continue by exhibiting an elliptic fibration on $X$
whose Jacobian visibly is a  K3 surface.

We start by considering the pencil of planes through the line $(z_3=z_4=0)$.
This yields a pencil of quartics with a singular point.
The pencil of lines through the singular point leads to the Weierstrass form
\begin{eqnarray*}
y^2 &  =  &
9x^6t^6-24x^6t^5+16x^6t^4-24x^5t^6+38x^5t^5-8x^5t^4+16x^4t^6-8x^4t^5-9x^4t^4-8x^4t^3\\
&& \;\;\;-16x^4t^2-8x^3t^4+20x^3t^3+8x^3t^2-16x^2t^4+8x^2t^3-9x^2t^2+32x^2t+32xt^2-58xt+9
\end{eqnarray*}
Note the symmetry in $x,t$.
Let $u=xt$.
Projection onto the projective $u$-line induces an elliptic fibration on $X$
with reducible fibers of type $I_6$ at $0, \infty$, $I_0^*$ at $-1$ and $I_2$ at $1$.
We did not find a section for this fibration, 
so we switch to its Jacobian $Jac(X)$
(which shares the same Hodge numbers,
 even the Picard number (over $\bar\QQ$) and the transcendental Galois representation over $\QQ$).
We find the following Weierstrass form for $Jac(X)$ with 2-torsion section $(0,0)$:
\[
Jac(X):\;\;
y^2 = 
x(x^2+
9(9u^2-58u+9)(u+1)^2x/4
+ 8100(u+1)^2u^2).
\]
By inspection of the degrees of the coefficients, $Jac(X)$ is a K3 surface.
The involution $u\mapsto 1/u$ of the base extends to $Jac(X)$ to define two involutions
connected by the hyperelliptic involution.
One has quotient a rational elliptic surface with reducible fibers of type $I_6, III$.
In consequence, the Mordell-Weil rank is two.
By the Shioda-Tate formula we deduce $\rho(Jac(X))\geq 19$ as required.

\subsection{2nd K3 surface}

The other involution on $Jac(X)$ is a Nikulin involution;
its quotient is an elliptic K3 surface $X'$ with singular fibers of type $I_6, I_1^*, III^*$ and a 2-torsion section.
In the new parameter $s=u+1/u$, it is given as
\begin{eqnarray}
\label{eq:X'}
X':\;\;
y^2 = 
x(x^2+(9s-58)(s+2)^2(s-2)x/4+8100(s+2)^3(s-2)^2).
\end{eqnarray}
We shall use that $X'$ fits into a Shioda--Inose structure;
more precisely $X'$ is sandwiched by a Kummer surface of product type
by means of Nikulin involutions.
To see this, it suffices to exhibit an elliptic fibration on $X'$ with section and two fibers of type $II^*$
(sometimes also called Inose pencil).
Abstractly this can be achieved by working out disjoint divisors of the given Kodaira type;
explicitly, however, the linear systems involved can become quite complicated, so we proceed in two steps by extracting one fiber of type $II^*$ at a time.

\subsection{Alternative elliptic fibration}

We start by extracting a singular fiber of type $IV^*$ from the above fibration as follows:
omit the opposite component of the $I_6$ fiber extended by zero section and identity component of the $I_1^*$ fiber.
The remaining fiber components disjoint from this and the 2-torsion section combine for divisors of type $\tilde E_8$ and $A_3$,
and  there still is a section.
In particular, we find 
\[
\NS(X') = U+A_3+E_6+E_8 \;\;\; \text{ and }\;\;\; T(X') = U+\langle 12\rangle
\]
(the latter through the discriminant form after Nikulin).
In terms of the above model \eqref{eq:X'},
the new fibration is given by the elliptic parameter $u=x/(s-2)$.
After rescaling $u=720/t$, a Weierstrass form is given by 
\[
X':\;\; y^2 = 
x^3
-216t^3(17t-15)x
-27(3375t^3-14393t^2+16965t-6075)t^4
\]
with reducible fibers $IV^*$ at $0$, $I_4$ at $1$ and $II^*$ at $\infty$.

\subsection{Inose pencil}
\label{ss:Inose}

We continue by extracting another divisor of Kodaira type $II^*$
from all but one far component of the $IV^*$ fiber extended by zero section and two components of the $I_4$ fiber.
Then the opposite component of the $I_4$ fiber serves as a section, 
and there is a disjoint divisor of type $E_8$ coming from the original $II^*$ fiber.
As elliptic parameter, one finds
\[
u = \frac{y+t(x-12t^2)(4-10t)+t^2(t-1)(1125t^2+225t-810)/2}{t^4(t-1)^2}.
\]
After some coordinate transformations over $\QQ$, this leads to the Weierstrass form of Inose's pencil
\begin{eqnarray}
\label{eq:Inose}
X':\;\;
y^2 = 
x^3
-436u^4x/3+u^5(5u^2/4-18997u/27-62500).
\end{eqnarray}
Following Inose \cite{Inose}, it is known 
that the quadratic base change $u=t^2$ leads to a Kummer surface of product type.
One can easily compute the j-invariants as
\[
j=16384/5,\;\;\; j'=-2^4109^3/5^6.
\]
Note that $j$ is the j-invariant of $E$ from \ref{ss:f},
and $j'$ describes an elliptic curve $E'$ which is 6-isogenous to $E$.
We shall now work out the above correspondence over $\QQ$.

\subsection{Kummer surface}

Recall that $E$ has a rational 6-torsion point.
The corresponding 6-isogeny leads to the minimal model  of $E'$:
\[
E':\;\; 
y'^2 = x'(x'^2+22x'+125).
\]
The Kummer surface $\Km(E\times E')$ can be given as cubic pencil
\[
\Km(E\times E'):\;\;\;
x(x^2+x-1) t^2 = x'(x'^2+22x'+125).
\]
Here projection onto the  $t$-line induces an elliptic fibration 
with fibers of type $IV^*$ at $0, \infty$.
We choose the base point $(0,0)$ of the cubics as zero section
and compute the Weierstrass form
\[
y^2 = 
x^3
-436t^4x/3+t^4(5t^4/4-18997t^2/27-62500).
\]
This is exactly the quadratic base change of the Inose pencil \eqref{eq:Inose} on $X'$
alluded to in \ref{ss:Inose}.
It follows that  all K3 surfaces $X, X'$ and $\Km(E\times E')$ share the same transcendental motive $V$ over $\QQ$.
Since $E$ and $E'$ are isogenous over $\QQ$ (and non-CM), 
we have in fact $V=\Sym H^1(E)$.
This concludes the proof of Theorem \ref{ss:thm}.
\qed

\subsection{Remark}

Note that the Tate conjecture is known for elliptic K3 surfaces with section.
In the present situation, we deduce that the Tate conjecture holds for $\tilde S$.
At a prime $p$ of good reduction, we find
\[
\rho(\tilde S\otimes\bar\FF_p) = 
\begin{cases}
51 & \text{ if $E$ has ordinary reduction at } p;\\
61 & \text{if $E$ has supersingular reduction at } p.
\end{cases}
\]
By work of Elkies, it follows that $\tilde S$ has infinitely many supersingular primes.
For instance, at $p=11$, $\tilde S$ attains supersingular reduction.
It would be interesting to see the additional divisors explicitly
and to have a moduli theoretic interpretation.

\section{Hilbert modular surfaces and Mumford-Tate groups}
\label{ss:HMS}

\subsection{The Shimura surface as Hilbert modular surface}
The Shimura surface is somewhat similar to a Hilbert modular surface 
described by Hirzebruch in \cite{Hi}: for example, both have fives cusps which are resolved by a cycle of three rational curves and the alternating group $A_5$
acts on both. However, the example from \cite{Hi} is a rational surface, whereas the Shimura surface is of general type. We will now show that the Shimura surface is dominated by a Hilbert modular surface, but it seems rather hard to determine the corresponding congruence subgroup explicitly.

The abelian surface $B_0=J(C_0)$ 
is the quotient of the vector space
$V:=\ZZ[\zeta]\otimes_\ZZ\RR$, endowed with a complex structure $J_0$,
by the lattice $\ZZ[\zeta]$ and the principal polarization on $B_0$ is given by
an alternating, $\ZZ$-bilinear form $E_0:\ZZ[\zeta]\times\ZZ[\zeta]\rightarrow\ZZ$. 

The complex structure $J_0$ on $V$, which defines $B_0$, can be deformed
to a complex structure $J$ in such a way that $E_0$ remains a principal polarisation (these deformations are parametrised by the Siegel space $\HH_2$). Moreover, we require that $J$ still commutes with multiplication by $\eta:=\zeta+\zeta^{-1}\in\ZZ[\zeta]$. 
These deformations are parametrised by a submanifold 
$\HH_1\times\HH_1\subset \HH_2$ which maps onto a Hilbert modular surface 
for the field $\QQ(\eta)\cong\QQ(\sqrt{5})$ in the moduli space $\cA_2$ of 
2-dimensional ppav's, cf.\ \cite{vdG}, Chapter IX. 
Let $B$ be the abelian surface defined by $J$, 
then $B$ is a principally polarised abelian surface and $\ZZ[\eta]\subset End(B)$. 
In particular, any matrix $M\in M_2(\ZZ[\eta])$ induces an endomorphism of $B^2$. 

We define an endomorphism $\psi_0$, in fact an isogeny, of 
$B_0^2$ by $\psi_0:(x,y)\mapsto (x+\delta y,x-\delta y)$
where $\delta=\zeta-\zeta^{-1}\in \ZZ[\zeta]=End(B_0)$.
This isogeny induces a homomorphism on the period lattice, 
given by the $2\times 2$ matrix $\Psi\in M_2(\ZZ[\zeta])$. 
One easily checks that there is a commutative diagram 
(note that $\delta^2=-3-\eta$):
{\renewcommand{\arraystretch}{1.3}
$$
\begin{array}{rcl}
\ZZ[\zeta]^2\,&\stackrel{\Psi}{\longrightarrow}&\,\ZZ[\zeta]^2\\
\stackrel{M}{\phantom{}} \downarrow\quad&&\quad\downarrow 
\stackrel{2\Phi_0}{\phantom{}}\\
\ZZ[\zeta]^2\,&\stackrel{\Psi}{\longrightarrow}&\,\ZZ[\zeta]^2,
\end{array}\qquad
\Psi\,=\,\begin{pmatrix} 1&\delta\\1&-\delta \end{pmatrix},\quad
M\,=\,\begin{pmatrix} \eta& -\eta-3\\ 1&\eta\end{pmatrix},\quad
\Phi_0\,=\,\begin{pmatrix} \zeta& 0\\ 0&\zeta^{-1}\end{pmatrix}.
$$
} 

The pull-back of the product polarisation $\Psi^*(E_0\oplus E_0)$ 
is obviously a polarisation on $B_0^2$.
We claim that it also defines a polarisation on $B^2$, for the deformations $J$ of the complex structure $J_0$ as above. 
To see this, we use the identification of $\NS(B_0^2)$ with the 
Rosati-invariant endomorphisms $End^s(B_0^2)\subset M_2(\ZZ[\zeta])$ 
(\cite{BL}, Prop.\ 5.2.1) given by $L\mapsto \phi_{L_0}^{-1}\phi_L$ 
where $L_0$ is a line bundle on $B_0^2$ with Chern class the principal polarisation $E_0\oplus E_0$. As $E_0\oplus E_0$ also defines a principal polarisation on $B^2$, 
we get an identification of $\NS(B^2)$
with $End^s(B^2)\subset M_2(\ZZ[\eta])\subset M_2(\ZZ[\zeta])$. 
The Rosati involution for an endomorphism defined by a matrix $\Omega\in M_2(\ZZ[\zeta])$ is given by 
$\Omega\mapsto \Omega'\,=\,{}^t\overline{\Omega}$ (the complex conjugation comes from \cite{BL}, Lemma 5.5.4, 
the transpose comes from pull-back of line bundles on $\widehat{B_0^2}=Pic^0(B_0^2)$, the dual of $B_0^2$). As we use $E_0\oplus E_0$ to define the Rosati involution, this polarisation corresponds to the identity in $End(B_0^2)$ 
and then $\Psi^*(E_0\oplus E_0)$ corresponds to
$\psi_0'\psi_0$ by \cite{BL}, Lemma 5.2.6. 
One easily computes that ${}^t\overline{\Psi}\Psi$ lies in $M_2(\ZZ[\eta])$, 
hence the polarisation $\Psi^*(E_0\oplus E_0)$ deforms to the class of a line bundle $L$ on $B^2$. It is still a polarisation, since this is an open condition (one now only needs $E(v,Jv)>0$).
 
Given a deformation $B$ of $B_0$ as above, we define an abelian variety $A$ as the image of $B^2$ under the isogeny determined by the matrix $\Psi$.
In particular, as $B$ is the quotient of the complex space $(V,J)$ by $\ZZ[\zeta]$, we can identify $A$ with the quotient of $(V,J)^2$ by the larger lattice $\Psi^{-1}(\ZZ[\zeta]^2)$.
As $\psi_0(B_0^2)=B_0^2$, the abelian variety $A$ is also a deformation of 
$B_0^2$. As $E_0\oplus E_0$ pulls back along $\Psi$ to a polarisation on $B^2$,
it defines a polarisation on $A$, which is principal since it is so on $B_0^2$.
The diagram shows that $A$ has an endomorphism induced by $2\Phi_0$. This endomorphism factors over the multiplication by two map on $A$,
$A\stackrel{2}{\longrightarrow} A\stackrel{\phi}{\rightarrow}A$,
and in this way it defines an endomorphism of $\phi$, 
which is of course induced by $\Phi_0$. 
In particular, $(A,\phi)$ is a deformation of $(B_0^2,\phi_0)$.

In this way we obtain a two dimensional (parametrised by $\HH_1\times \HH_1$)
family of abelian fourfolds with automorphism of order five $(A,\phi)$ which are deformations of $(A_0,\phi_0)$. 
On the other hand, the points in the Shimura surface $S\subset \overline{\cA_4(2,4)}$ also
parametrise principally polarised abelian fourfolds $A$ with an automorphism 
$\phi$ which are deformations of $(A_0,\phi_0)$. 
These surfaces therefore coincide and thus $S$
is dominated by a Hilbert modular surface belonging to the field $\QQ(\zeta)$. 
The precise subgroup $\Gamma$ of $\mbox{SL}_2(\QQ(\sqrt{5}))$ such that $S$ is the compactification of
$\Gamma\backslash(\HH_1\times\HH_1)$ can, in principle, 
be determined from the data given here.

\subsection{Mumford-Tate groups}
We briefly review the Mumford-Tate groups of the abelian varieties considered here.
This provides another way to see that the Shimura surface is dominated by a Hilbert modular surface and, in principle, it allows one to determine the cocompact discrete subgroup $\Gamma$ of $SL(2,\RR)$ such that $\Gamma\backslash \HH_1$ is isomorphic to the Shimura curve.

Let $(X,L,\phi)$ be a four dimensional ppav with an automorphism $\phi$ of order five as in section \ref{geno}. 
Then $X\cong (\Lambda\otimes\RR)/\Lambda$ where $\Lambda\cong \ZZ[\zeta]^2$
and the action of $\phi$ on $\Lambda\cong \ZZ[\zeta]^2$  is given by multiplication by $\zeta$. The line bundle $L$ defines a bilinear alternating form $E$ on $\Lambda$,
its first Chern class, which is given by a skew Hermitian $2\times 2$ matrix $T$
with coefficient in $K:=\QQ(\zeta)$ (cf.\ \cite{BL}, Ch.\ 9):
$$
E\,:\,\Lambda\,\times\,\Lambda\,\longrightarrow\,\ZZ,\qquad
E(x,y)\,=\,tr_{K/\QQ}({}^txT\overline{y}),\qquad (x,y\in \Lambda\cong \ZZ[\zeta]^2).
$$
The complex structure $J$ on $V_\RR:=\Lambda\otimes\RR$ which defines $X$ is symplectic ($E(Jx,Jy)=E(x,y)$) and commutes with the action of $K$.

The special Mumford Tate group of the rational polarised Hodge structure
$(V:=\Lambda\otimes_\ZZ\QQ,\,E)$ is the smallest algebraic subgroup $G$ of $Sp(E)$, defined over $\QQ$, such that $a+bJ\in G(\RR)$ for all $a,b\in \RR$ with $a^2+b^2=1$. In what follows we will often `confuse' an algebraic group defined over $\QQ$ with the group of its $\QQ$-points.

As the actions of $J$ and $K$ on $V_\RR$ commute, $a+bJ$ lies in $L(\RR)$ where $L$ is the algebraic subgroup of $Sp(E)$ of elements which commute with the $K$-action on $V$ (this group is obviously defined over $\QQ$). In particular, $G\subset L$. The group $L$ can be shown to be isomorphic to a unitary group:
$$
L\,:=\,\{A\in GL(V)\,:\, E(Ax,Ay)\,=\,E(x,y),\quad 
aA=Aa\quad \forall a\in K\,\}
\,=\,U(H_K),
$$
where we view $V$ as a two dimensional $K$-vector space and 
the hermitian form $H$ is obtained from the
skew-Hermitian matrix $T$:
$$
H\,:\,V\times V\,\longrightarrow \, K,\qquad
H(x,y)\,:=\, \delta{}^txT\overline{y},\qquad \delta:=\zeta-\overline{\zeta}_5.
$$
and $U(H)$ is the group of $K$-linear maps on $V$ which preserve $H$.

The special Mumford-Tate group $G$ is contained in $SU(H)$, where
$$
SU(H)\,=\,\{A\,\in\, U(H)\,:\,\det A\,=\,1\}.
$$  
Let $K_0:=\QQ(\sqrt{5})$.  The isomorphism 
$$
K_0\otimes_\QQ\RR\,\stackrel{\cong}{\longrightarrow}\, \RR\times\RR,\qquad
\sqrt{5}\otimes 1\,\longmapsto\,(\sqrt{5},-\sqrt{5})
$$
induces an isomorphism $V_\RR\cong V_1\oplus V_2$, with each $V_j\cong\RR^4$, 
and an isomorphism
of real Lie groups
$$ 
SU(H)(\RR)\,\cong\,\,SU(r_1,s_1)\,\times SU(r_2,s_2),\qquad 
(r_1+s_1=r_2+s_2=2).
$$
The $r_j$ are determined as follows: the action of $K\otimes\RR\cong \CC\times\CC$ 
on $V=V_1\oplus V_2$ defines the structure of a complex vector space on the $V_j$. The complex structure $J$ is $\CC$-linear and $r_j$ is the dimension of the 
eigenspace where $J=+i$ in $V_j$ (note that $SU(r,s)\cong SU(s,r)$). 

The corresponding Shimura variety is of the form $\Gamma\backslash SU(H)(\RR)/K(SU(H))$
where $K(SU(H))$ is a maximally compact subgroup and $\Gamma$ is a discrete subgroup.
For the Shimura curve we have $(r_1,s_1)=(2,0)$ (note that $SU(2,0)$ is compact) and $(r_2,s_2)=(1,1)$ whereas for the Shimura surface we have $(r_1,s_1)=(r_2,s_2)=(1,1)$. As $SU(1,1)\cong SL(2,\RR)$ 
and $K(SU(1,1))\cong SO(2,\RR)$ we find that the the $d$-dimensional Shimura variety is 
a quotient of $\Gamma_d\backslash\HH_1^d$ for $d=1,2$ where $\Gamma_d$ is a discrete subgroup of $SL(2,\RR)^d$.

To describe $\Gamma_d$ in more detail, we choose a $K$-basis of $V$ 
such that $T\,=\,\mbox{diag}(\delta_1,\delta_2)$ for certain $\delta_j\in K$.

The unitary group of the Hermitean form defined by $H:=\delta T$
in $GL(2,K)$ consists of the matrices $A$ such that
$$
{}^t\overline{A}(\delta T)A\,=\,\delta T\quad\Longleftrightarrow\quad
{}^t\overline{A}TA\,=\,T\quad\Longleftrightarrow\quad
A^{-1}\,=\,T^{-1}({}^t\overline{A})T.
$$
For any invertible $2\times 2$ matrix $A$ with $\det A=1$ one has
$$
A^{-1}\,=\,S({}^tA)S^{-1},\qquad\mbox{with}\quad S\,=\,
\begin{pmatrix} 0&1\\-1&0\end{pmatrix}.
$$
Thus for $A$ in the special unitary group we have
$$
A\,=\,U\overline{A}U^{-1},\qquad 
U\,=\,{}^t(TS)\,=\,\begin{pmatrix} 0&-\delta_2\\ \delta_1&0\end{pmatrix},
$$
or equivalently, $AU=U\overline{A}$.
Thus we  consider the set
$$
D\,:=\,\{\,A\in\,M_2(K)\,:\,A=U\overline{A}U^{-1}\,\}.
$$
Note that $D$ is a linear space over the real subfield $K_0$ of $K$
and that $A,B\in D$ implies $AB\in D$, hence $D$ is a $K_0$-algebra.
A $K_0$-basis of $D$ is easy to find explicitly from the condition $AU=U\overline{A}$:
$$
D=K_0 1\oplus K_0{\bf i}\oplus K_0{\bf j}\oplus K_0{\bf k},\qquad
{\bf i}:=\begin{pmatrix} \delta_1 &0\\0&-\delta_1\end{pmatrix},\quad
{\bf j}:=\begin{pmatrix} 0 &\delta_2 \\-\delta_1&0\end{pmatrix},
\quad {\bf k}:={\bf i}{\bf j}.
$$
Multiplication in $D$ satisfies:
$$
{\bf i}^2\,=\,\delta_1^2,\quad{\bf j}^2\,=\,-\delta_1\delta_2,\quad
{\bf i}{\bf j}=-{\bf j}{\bf i}
$$
note that $\delta_1^2,-\delta_1\delta_2\in K_0$. 
Hence $D$ is the quaternion algebra over $K_0$ which is usually denoted
by
$$
D\,=\,(\delta_1^2,-\delta_1\delta_2)_{K_0}.
$$
The quaternion algebra $D$ has a canonical (anti)involution defined by
$$
\overline{\alpha}\,=\,
\overline{a_0+a_1{\bf i}+a_2{\bf j}+a_3{\bf k}}\,:=\,
a_0-a_1{\bf i}-a_2{\bf j}-a_3{\bf k}
$$
and one verifies that
$$
N(\alpha)\,:=\,\alpha\overline{\alpha}\,=\,
a_0^2-\delta^2a_1^2+\delta\delta'a_2^2+\delta^3\delta'a_3^2,
$$
thus $\alpha$ is invertible in the algebra $D$ iff $N(\alpha)\neq 0$.
In particular, we have a homomorphism of multiplicative groups
$$
N\,:\,D^\times\,\longrightarrow\,K_0^\times,\qquad 
\alpha\,\longmapsto\,N(\alpha)\,:=\,\alpha\overline{\alpha}
$$
It is easy to check that if $\alpha$ is given by the matrix $A\in M_2(K)$ then $\overline{\alpha}$ is given by ${}^t\overline{A}$ and that
$N(\alpha)=\det(A)$.
Moreover, if $A\in D^\times$ and $\det A=1$ then, using $A^{-1}=S({}^tA)S^{-1}$, one obtains that ${}^t\overline{A}TA=T$, that is, $A$ is in the unitary group of $H_K$, hence we have an isomorphism of algebraic groups $SU(H)\cong D^\times_1$.

Thus the Shimura curve, $\Gamma_1\backslash\HH_1$, is defined by a subgroup $\Gamma_1\subset D_1^\times$. It is isomorphic to $\PP^1$ and thus has genus zero. In \cite{J} genus zero quotients $\Gamma\backslash\HH_1$, where $\Gamma$ is the group of units of a maximal order in a quaternion algebra with center a real quadratic field are classified. Note however that it is not clear that $\Gamma_1$ is of this type, probably it has finite index in such a group.

In case of the Shimura surface, one can check that the quaternion algebra is actually isomorphic to the matrix algebra  $M_2(K_0)$ and thus $\Gamma_2$ is subgroup of $SL_2(K_0)$, which implies that the Shimura surface is dominated by a Hilbert modular surface. 

We observe that explicit equations for Hilbert modular surfaces were recently investigated in \cite{Carls}. Explicit genus two curves whose Jacobians have endomorphisms by $\ZZ[\sqrt{5}]$ were studied in \cite{W}.

\subsection{Remark}

Our interpretation of $\tilde S$ as an Hilbert modular surface 
agrees very well with our results on the zeta function.
Namely work of Oda \cite{O} asserts that to a primitive Hilbert modular form of weight 2 with respect to $\Gamma_0(\mathfrak n)$,
one can associate a Hodge structure $M_f$ of weight two
and elliptic curves $E, E'$ over $\QQ$
such that there is an isomorphism of Hodge structures $M_f \cong H^1(E)\otimes H^1(E')$.

Note that even if this result should apply to the present situation,
the Galois representations would a priori only be isomorphic over some finite extension of $\QQ$.
That is the reason why we still need the explicit geometric relation over $\QQ$
between $\tilde S$ and a Kummer surface of product type in order to prove Theorem \ref{ss:thm}.

\subsection*{Acknowledgements}

We thank Fabrizio Catanese and Jan-Christian Rohde for helpful discussions.
This project was started when the second author held a position
at the University of Copenhagen.
Part of this research was conducted when the second author visited
the Dipartimento di Matematica of Universit\`a di Milano.
He is very grateful for the kind hospitality.

\end{document}